\documentclass{elsart}
\usepackage{ifpdf}
\usepackage{graphicx,amssymb,amsmath}%
\usepackage[square]{natbib}

\ifpdf
\usepackage[
pdfstartview =FitH,
a4paper=false,
  breaklinks=true,
  colorlinks=true,%
  linkcolor=blue,anchorcolor=blue,%
  citecolor=blue,
  menucolor=blue,
  urlcolor=blue]{hyperref}
\else
\fi

\usepackage[infoshow]{tracefnt}
\def\R{{\mathbb{R} }}
\newcommand{\RE}{\operatorname{\mathcal{R}e}}

\theoremstyle{definition}

\newtheorem{theorem}{Theorem}[section]
\theoremstyle{remark}
\newtheorem{remark}{Remarks}[section]

\newenvironment{proof}[1][Proof]{\begin{trivlist}
\item[\hskip \labelsep {\bfseries #1} ]}
 {\end{trivlist}}
\def\1{1\kern-.20em {\rm l}}

\def\EE{\mathbb{E}}
\def\S{\mathbb{S}}

\let \a=\alpha
\let \b=\beta
\let \G=\Gamma
\let \g=\gamma
\let \l=\lambda
\let \s=\sigma
\let \dy=\displaystyle

\let \t=\tau
\let \th=\theta

\renewcommand{\d }{\delta}
\let \text=\hbox

\makeatletter
\def\elsartstyle{%
    \def\normalsize{\@setfontsize\normalsize\@xiipt{14.5}}
    \def\small{\@setfontsize\small\@xipt{13.6}}
    \let\footnotesize=\small
    \def\large{\@setfontsize\large\@xivpt{18}}
    \def\Large{\@setfontsize\Large\@xviipt{22}}
    \skip\@mpfootins = 18\p@ \@plus 2\p@
    \normalsize
}
\@ifundefined{square}{}{\let\Box\square}
\makeatother
\journal{  }

\bibpunct{[}{]}{,}{n}{}{;} 
\begin{document}

\begin{frontmatter}
\title{Spectral representation of some non stationary $\a$-stable processes}

\author{Nourddine Azzaoui}
\address{Institut de Mathématiques de Bourgogne (IMB)\\
9 avenue Alain Savary\\
21000 Dijon, France}
\date{ }
\ead{nourddine.azzaoui@math-univ-lille1.fr}

\begin{abstract}
In this paper, we give a new covariation spectral representation of some non stationary symmetric $\alpha$-stable processes (S$\a$S). This representation is based on a weaker covariation pseudo additivity condition which is more general than the condition of independence. This work can be seen as a generalization of the covariation spectral representation of processes expressed as stochastic integrals with respect to independent increments S$\a$S processes (see Cambanis (1983)) or with respect to  the general concept of independently scattered S$\a$S measures (Samorodnitsky and Taqqu 1994). Relying on this result we investigate the non stationarity structure of some harmonisable S$\a$S processes especially those having periodic or almost-periodic covariation functions.
\end{abstract}

\end{frontmatter}

\section{Introduction}
In this work, we are interested in a family of stochastic processes having infinite second order moments. It is about the class of symmetric $\a$-stables (S$\a$S) processes, with ($1< \alpha<2$), having the stochastic integral representation:
\begin{equation}\label{eq2-1}
X_t=\int_{-\infty}^{\infty}{f(t,\lambda)d\xi(\lambda)},
\end{equation}
where $\xi$ is an S$\a$S stochastic process. By using the concept of the covariation that will be recalled later in (\ref{eq25}) or (\ref{eq26}) and its relationship  with the convergence in probability, \citet{Cam83} have given general conditions for the existence of the integral (\ref{eq2-1}), they are: 
\begin{itemize}
\item {\bf (I)} The process $\xi$ is right continuous with respect to the convergence in probability.
 \item {\bf(II)} for all linear combination $\zeta$ of increments of $\xi$, the map $v : t\longmapsto[\xi(t),\zeta]_{\alpha}$ is of bounded variations where $[. , .]_{\alpha}$ denote the covariation.
\end{itemize}
Due to their impulsive nature, these processes provide appropriate models in various application fields: for example in financial and econometric modeling see \citet{Uch99}, communications, signal processing see \citet{Nik96}, nuclear physics, astronomy see \citet{Uch99} and references within... As an alternative to the covariance, which is undefined in this case, the covariation was introduced by \citet{Mil77}. Although this new dependence measure was conceived to substitute the covariance when $1<\a<2$, it is not as convenient because it does not have some of its suitable properties. One of its great defects is that it is not, in general, additive with respect to its second variable. This last property play an important role in integral representation of the covariation with respect to a spectral measure (control measure). In the case of symmetric  $\a$-stables processes, \citet{Cam83} have used the concept of independence as a sufficient condition  for the covariation additivity what have enabled him to build a spectral representation of the covariation. This case have been widely studied in literature and was the keystone of several works in theoretical as applied fields. Among others, it was used in linear regression (\citet{Sam94}), prediction problems (\citet{Cam89},\cite{Cam84})and spectral analysis of stationary processes see \citet{Sam94} and reference within.

In this paper, we generalize works of \citet{Cam83} by replacing the increments independence condition on the process $\xi$, by a weaker one permitting the additivity of the covariation. Using this property we give, in the general sense of \citet{Mor56}, a new integral representation of the covariation function of the process $X$ with respect to a characterizing bimeasure. This representation which is similar to the covariance spectral representation of second order processes, will play an important role in the study of a wide class of non-stationary $\a$-stables processes. For instance, the class of second order harmonisable processes, introduced by Loeve, have been extensively studied by \citet{Rao82, Rao86}. In particular, periodically and almost periodically correlated second order processes have been thoroughly investigated, see for instance \citet{Hur89,Hur91, Deh94, Deh93}. In the Cambanis's spectral representation, where $\xi$ have independent increments, all harmonisables processes are covariation stationary: their covariation depend only on time difference, see \citet{Cam83}. However, our spectral representation give an understanding of the exact manner in which the increments of the processes $\xi$ are not independents gives some information about the non-stationarity nature of the process $X$. In particular, similarly to second order processes, we study an important class of  harmonisables non stationary processes for which the covariation is periodic or almost periodic that we call (\textit {periodically and almost-periodically covariated processes}). We show that some fundamental results concerning second order periodically and almost periodically correlated processes remain true for our periodically and almost-periodically covariated $\a$-stable processes. This paper is organized as follows: we begin by a brief overview on stable variables and processes and some preliminary results. We give, in theorem \ref{th2-1}, our first result concerning a weaker condition for covariation additivity. After a brief recall on Morse and Transue integration,  we provide our main spectral representation. finally, a classification of harmonisable processes according to the structure of their characterizing bimesure is provided. All the proofs will be given in appendix.

\section{Definitions, notations and preliminary results}
 A real centered random vector $X^d$=($X_1,X_2,\dots,X_d$)  is  symmetric $\alpha$-stable (S$\alpha$S) if and only if its characteristic function is given by: 
\begin{equation}\label{eq-car-real}
\phi_{X^d}(t_1,\dots,t_d)=\exp\{-\int_{\S_d}{\left|\sum_{i=1}^{d}{t_i 
s_i}\right|^{\alpha}d\Gamma_{(X_1,\dots,X_d)}(s_1,\dots,s_d)}\},
\end{equation}
where $\Gamma_{(X_1,\dots,X_d)}$ is an unique finite symmetrical measure defined on the unit sphere $\S_d$ of $\R^d$, see \citet{Sam94}. When $X^d$ is a complex S$\alpha$S vector, its characteristic function is defined, for all complex $\theta=(\theta_1, \dots, \theta_d)$ in $\mathbb{C}^d$, by: 
\begin{equation}\label{eq2-2}
\Phi_X(\theta_1, \dots, \theta_d)= \mathbb{E}\left(\exp\{\imath
\RE\langle
 \theta , \overline{X} \rangle\}\right) =\mathbb{E}\left(\exp\{\imath \RE\sum_{k=1}^d{\theta_k \overline{X}_k}\}\right).
\end{equation}
The characteristic function (\ref{eq2-2}) is computed like (\ref{eq-car-real}) applied to the real S$\a$S random vector $X^{2d}=(X^1_1,X^2_1,\dots,X^1_d, X^2_d)$ where $X_j=X^1_j+ \imath X^2_j$. In this paper, when $X^{^d}$ is complex we denote $\G_{X^{d}}$ the unique spectral measure corresponding to the real vector $X^{2d}=(X^1_1, X^2_1,\dots,X^1_d, X^2_d)$. Let $\phi$ be the Fourier transform of the spectral measure $\G_{X^d}$. When $X^d$ is real, it is given for all $\theta=(\theta_1,\dots,\theta_d) \in \mathbb{R}^d$ by:
\begin{equation}\label{eq2-17}
\phi(\theta)=\int_{\S_{d}}{\cos\left(\sum_{i=1}^d{\theta_is_i}\right)d\Gamma_{_{X^d}}(s_1,\dots,s_d)}.
\end{equation}
For the complex case, it is given by,
\begin{equation}\label{eq2-17b}
\phi(\theta)=\int_{\S_{2d}}{\cos\left(\RE\sum_{i=1}^d{\theta_i
\overline{s}_i}\right)d\Gamma_{_{X^d}}},
\end{equation}
where $s_j=s_j^1+\imath s_j^2$ and $\theta=(\theta_1,\dots,\theta_d)\in \mathbb{C}^d$. 
\begin{remark}
The finite measure $\G_{X^d}$ is defined on the unit sphere which is compact this implies that it have finite moments of all orders, by the sequel $\phi$ is, in particular  three times differentiable and we have:
\begin{equation}\label{eq2-18}
\frac{\partial^3\phi}{\overline{\partial}\theta_i\overline{\partial}\theta_j\overline{\partial}\theta_k}(\theta_1,...,\theta_d)=\int_{\S_{d}}{s_i s_j s_k\sin\left(\RE\left(\sum_{\ell=1}^d{\theta_{\ell} s_{\ell}}\right)\right)d\Gamma_{_{X^d}}(s_1,\dots,s_d)},
\end{equation}
where $\dy \frac{\partial}{\overline{\partial}\theta_j}=
\left(\frac{\partial}{\partial\theta_j^1} + \imath
\frac{\partial}{\partial\theta_j^2}\right)$ and $\dy \frac{\partial}{{\partial}\theta_j}=
\left(\frac{\partial}{\partial\theta_j^1} - \imath
\frac{\partial}{\partial\theta_j^2}\right)$ for $\theta_j =
\theta_j^1 + \imath \theta_j^2$. When $\theta_j$ is real then these operators are the usual partial derivatives.
\end{remark}
The covariation was conceived to replace the covariance which is undefined for S$\a$S variables. It is given in the real case by:
 \begin{equation}\label{eq25}
[X_1,X_2]_{\a}=\int_{\S_2}{s_1 s_2^{<\a-1>}d\Gamma_{(X_1,X_2)}(s_1,s_2)},
\end{equation}
where $s^{<\b>}=\mathrm{sign}(s)  |s|^{\b}$, see \citet{Sam94}. Similarly for $X_1=X^1_1 + \imath X_1^2$ and
$X_2=X_2^1 + \imath X^2_2$, the covariation of $X_1$ on $X_2$ is given by:
\begin{equation}\label{eq26}
[X_1,X_2]_{\a}=\int_{\S_4}{(s^1_1+\imath s_1^2) (s_2^1+ \imath
s^2_2)^{<\a-1>}d\Gamma_{(X_1,X_2)}(s^1_1,s_1^2,s_2^1,s^2_2)},
\end{equation}
and $z^{<\b>}=|z|^{\b-1} \overline{z}$, see \citet{Cam83}. The last author have shown that integrals in (\ref{eq25}) and (\ref{eq26}) remain unchanged if we replace $\Gamma_{(X_1,X_2)}$ by a higher order $\Gamma_{(X_1,\dots,X_d)}$. This property is easily derived from the fact that, for all complexes ($a_j$) and ($b_j$) we have:
\begin{equation}\label{eq2-3}
\left[\sum_{j=1}^d{a_j X_j},\sum_{j=1}^d{b_j X_j}\right]_{\alpha}=\int_{\S_{2d}}{(\sum_{j=1}^d{a_j(s_j^1+\imath
s_j^2)})(\sum_{j=1}^d{b_j(s_j^1+\imath
s_j^2)})^{^{<\alpha-1>}}d\Gamma_{X^{d}}}.
\end{equation}
It was also shown that, for $1<\alpha\leq 2$, the map $X\longmapsto
\|X\|_{\alpha}=([X,X]_{\alpha})^{\frac{1}{\alpha}}$ is a norm on the vector space induced by the set of S$\alpha$S random variables. For more details and further properties see \citet{Cam83}.
\subsection{Sufficient condition for the additivity of the covariation}
 In this section we are interested in the additivity of the covariation function with respect to its second variable, that is $[Y,a X_i+ b X_j]_\alpha= [Y,a X_i]_\alpha+ [Y,b X_j]_\alpha$, for all $a$ and $b$. This property is usually provided by the independence of $X_i$ and $X_j$ or more generally by: $\Gamma\left\{(s_1,...,s_d)\in \S_d,
s_i.s_j\neq 0\right\}=0$, see \citet{Cam80, Cam81}. The idea of this paper is to introduce a new condition on the vector $ X^d$ more general than the condition of independence and permitting the additivity of the covariation. The next theorem provides an answer of this issue.
 
\begin{theorem}\label{th2-1}
For the covariation to be additive with respect to its second variable  that is, for all $i_0\in \{1,\dots,d\}$:
\begin{equation}\label{eq2-19}
\forall \theta_1,\dots, \th_d \in \mathbb{C}, \hspace{0.2cm}
[X_{i_0},\th_1X_1+\dots+\th_d X_d]_{_{\alpha}} =
[X_{i_0},\th_1X_1]_{_{\a}}+\dots+[X_{i_0},\th_d X_d]_{_{\a}}.
\end{equation}
It is sufficient that, for all $i, j$ and $k \in \{1,\dots,d\}$ not all equals, the Fourier transform $\phi$ fulfill the next condition:
\begin{equation}\label{eq2-20}
\forall \th_1,\dots, \th_d \in \mathbb{C}, \hspace{1.5cm}
\frac{\partial^3
\phi}{\overline{\partial}\th_i\overline{\partial}\th_j\overline{\partial}\th_k}(\th_1,\dots,\th_d)=0.
\end{equation}
\end{theorem}

{\bf Examples} 
\begin{itemize}
\item Let us remark that the result (\ref{eq2-20}) of theorem \ref{th2-1}
is a generalization of the independence condition. Indeed,
suppose that $X_1,\dots,X_d$ are pairwise independent real S$\a$S variables, then according to  [\citet{Sam94}, p. 68], the measure $\G_{X^d}$ is concentrated on the intersection points of the sphere $\S_d$ and the base axis of $\mathbb{R}^d$. We can then write $\G_{X^d}$ as:
$$
\G_{_{X^d}}=a_1 [\d_{_{(1,0,\dots,0)}}+\d_{_{(-1,0,\dots,0)}}]+\dots +a_d [\d_{_{(0,0,\dots,1)}}+\d_{_{(0,0,\dots,-1)}}],
$$
where $a_1,\dots,a_d$ are positifs weights. In this case, it is easy to see that the Fourier transform of  $\G_{X^d}$  can be given by:

$$
\phi(\th_1,\dots,\th_d)= \sum_{i=1}^{d}{a_i \cos(\th_i)}.
$$
If we partially differentiate this formula two times we find a null function. This implies clearly that the function $\phi$ given in this example satisfies the condition (\ref{eq2-20}).
\item \underline{Non trivial example}: To get a simple example of S$\a$S vector $X^d$ with spectral measure $\G_{_{X^d}}$ determined by its  Fourier transform $\phi$ and verifying the condition (\ref{eq2-20}), we choose $\phi$ of the form:
$$
\phi(\th_1,\dots,\th_d)=\sum_{i=1}^{d}{\varphi_{i}(\th_i)},
$$
where $\varphi_{i}$'s are three times differentiable even real functions. By applying the Bochner's theorem, in order that $\phi$ be a Fourier transform of a finite measure it is necessary and sufficient that it is positive definite. Then, for instance, one can take the $\varphi_{i,j}$'s as three times differentiable characteristic functions of random variables. 
\end{itemize}
\subsection{Bimeasure construction and Morse Transue integral}
We recall that a stochastic process $\xi=(\xi_t, t\in \R)$ is S$\a$S if all finite subset of $\xi$ is also an  S$\a$S vector. This is equivalent to the fact that all finite linear combinations of elements of $\xi$ are also S$\a$S random variables, see  \citet{Cam83, Sam94}.
Let us consider the increments application, mapping each interval $[s,t[$ to the complex  S$\a$S random variable  $d\xi$ defined by:
\begin{equation}\label{eq2-34}
d\xi\left([s,t[\right) = \xi_t-\xi_s.
\end{equation}
According to \citet{Cam83}, under the conditions (I) and (II), the application $d\xi$ may be extended to a complex random measure on the Borel $\s$-algebra $\mathcal{B}(\mathbb{R})$. For the construction of our bimeasure we will need the following condition.

\textbf{\underline{Condition $\mathcal{A}$}}\footnote{The letter $\mathcal{A}$ is for "Additivity" } :\\
{ \it We will say that $d\xi$ satisfy the condition $\mathcal{A}$ if and only if for all $n\geq 2$, and for all, pairwise distinct Borel sets $A_1,...,A_n$, the S$\a$S random vector $\dy
(d\xi(A_1),...,d\xi(A_n))$ satisfy the additivity condition (\ref{eq2-20}) of theorem \ref{th2-1}.}

Suppose now that the random measure $d\xi$ satisfy the condition $\mathcal{A}$
and consider the complex valued set function $F$ defined on
$\mathcal{B}(\mathbb{R})\times \mathcal{B}(\mathbb{R})$ by:
\begin{equation}\label{eq2-35}
 F :\begin{array}{lll}\mathcal{B}(\mathbb{R})\times \mathcal{B}(\mathbb{R}) & \longrightarrow & \dy \mathbb{C}\\
 \dy (A,B)& \longmapsto & \dy [d\xi(A), d\xi(B)]_{\alpha}
  \end{array}.
\end{equation}
The application $F$ is additive with respect to its two variables: it is a bimeasure. Indeed, 
the additivity of the first component comes from the linearity of the covariation with respect to its first variable see \citet{Sam94}. For the second component, let $B_1$ and $B_2$ be two disjoint Borel sets then for all fixed $A$ we have,
\begin{equation}\label{eq2-36}
F(A, B_1 \cup B_2)=[d\xi(A), d\xi(B_1 \cup B_2)]_{\alpha}=
[d\xi(A), d\xi(B_1) + d\xi(B_2)]_{\alpha}.
\end{equation}
Since $d\xi$ satisfy the condition $\mathcal{A}$ then,
$$
F(A, B_1 \cup B_2)=[d\xi(A), d\xi(B_1)
]_{\alpha}+[d\xi(A), d\xi(B_2)]_{\alpha} =F(A, B_1)+F(A,B_2).
$$
We can show that $F$ is also $\s$-additive with respect to its two variables. Indeed, according to
\citet{Rao82} or \citet{Rao86}, it suffices to show that, for all Borel sets $A_n$ decreasing to $\emptyset$, 
$F(A_n,A_n)$ converges to 0. Since $\xi$ is right continuous with respect to the weak topology of convergence in probability\footnote{ The topology induced by convergence in probability is equivalent to the covariation norm topology, see [\citet{Sam94}, p. 95].}, we deduce that if $(A_n)_n$ decreases to $\emptyset$ then $d\xi(A_n)$
converges in probability to 0. This implies that $F(A_n,A_n)=\|d\xi(A_n)\|_{\a}^{\a}$ converges to 0.

In the case of second order processes, it is known that the bimeasure associated with the covariance function is positive definite. This property plays an important role in the construction of integrals with respect to a bimeasure see \citet{Rao86}. In our case of  S$\a$S processes, the bimeasure $F$ defined in (\ref{eq2-35}) have the following property:
\begin{prop}
The bimeasure $F$ defined in (\ref{eq2-35}) have a property similar to Bochner's positive definiteness: for all complex $z_1, \dots, z_n $ and for all pairwise distinct Borel sets $ A_1, \dots, A_n $, we have:
\begin{equation}\label{eq2-bim}
\sum_{i=1}^n\sum_{j=1}^n{z_i (z_j)^{<\a-1>} F(A_i, A_j)} \geq 0.
\end{equation}
\end{prop}
The proof of this property is easy, it suffices to use the condition $\mathcal{A}$. Indeed,
 $$
 \begin{array}{lll}
 \dy \sum_{i=1}^n\sum_{j=1}^n{z_i (z_j)^{<\a-1>} F(A_i, A_j)}&= &\dy \sum_{i=1}^n\sum_{j=1}^n{z_i (z_j)^{<\a-1>} [d\xi(A_i), d\xi(A_j)]_{\a}},\\
 & = & \dy \left[\sum_{i=1}^n{z_i d\xi(A_i)}, \sum_{i=1}^n{z_i d\xi(A_i)} \right]_{\a},\\
 & = &\dy  \left\|\sum_{i=1}^n{z_i d\xi(A_i)}\right\|^{\a}_{\a} \geq 0.
 \end{array}
 $$
\par
In all the rest of the paper we suppose that the bimeasure $F$ satisfy a condition similar to bounded Fr\^echet variations but weaker than Vitali's\footnote{ In this paper we will not be interested in some complicated considerations, about different bounded variations nor in their relationship with the existence of different type of integration theory with respect to a bimeasure $F$, discussed in \citet{Rao86}.  We just recall that, $F$ have bounded Vitali variations if it verifies: 
\begin{equation}\label{eq2-37}
\|F\|(\mathbb{R}\times \mathbb{R})\triangleq \sup\left\{
\sum_{i=1}^{n} \sum_{j=1}^{n}{ |F(A_i,B_j)}|,  A_i\times B_j \hbox{ distincts } \right\}< \infty.
\end{equation}
}. It is given by:
\begin{equation}\label{Frech}
\sup\left\{
\sum_{i=1}^{n} \sum_{j=1}^{n}{a_i (a_j)^{<\a-1>} F(A_i,A_j)},
A_i\in \mathcal{B}(\mathbb{R}) \hbox{ and } |a_i|<1\right\}< \infty.
\end{equation}
Let $B$ be a fixed Borel set and consider the complex measure defined by, $F_B : A\longmapsto
F_B(A)=F(A,B)$. According to (\ref{Frech}) the complex measure $F_B$ have a finite total variation  on $\mathcal{B}(\mathbb{R})$ (for the definition of total variation see [\citet{Dun58}, p. 97], we will recall it later in (\ref{eq2-39})). As in \citet{Rao82}, we define, in the sense of \citet{Dun58}, the integral of a bounded function $f$ with respect to $F_B$ that we denote by:
\begin{equation}
  \tilde{I}_1(f , B)\triangleq \int_{\mathbb{R}}{f(\lambda)F_B(d\lambda)}.
 \end{equation}
Now, let $f$ be a fixed bounded map and consider the set function $B\longmapsto \tilde{I}_1(f , B)$ which is also a complex measure with bounded total variation on $\mathcal{B}(\mathbb{R})$. Similarly, we define the integral of a bounded function $g$
with respect to $\tilde{I}_1(f , .)$:
\begin{displaymath}
  I_1(f , g)\triangleq\int_{\mathbb{R}}{g(\lambda)\tilde{I}_1(f , d\lambda)}.
 \end{displaymath}
In the same way, by integrating first with respect to the complex measure 
$\tilde{F}_A(B)=F(A,B)$ and $A$ fixed, one can construct another integral $\dy I_2(f ,
g)=\int_{\mathbb{R}}{f(\lambda)\tilde{I}_2(d\lambda,g)}$ where $\dy
\tilde{I}_2(A,g)=\int_{\mathbb{R}}{g(\lambda)\tilde{F}_A(d\lambda)}$.
Note that in general the two integrals $I_1$ and $I_2$ are not equal (see counter example in \citet{Mor56}). When $I_1(f , g)= I_2(f,g)$ the couple $(f,g)$ is said Morse Transue integrable (MT-integrable) and the common value will be denoted as a double integral:
\begin{equation}\label{eq2-38}
I(f,g)=\int_{\mathbb{R}}\int_{\mathbb{R}}{f(\lambda)g(\lambda')F(d\lambda,d\lambda')}.
\end{equation}
As it was mentioned by \citet{Rao82}, it is not easy to identify the class of all  MT-integrable functions. It should be noted that complexes measures on $\mathcal{B}(\R^2)$ are particular bimeasures, in this case MT-integration is the same as the usual Lebesgue (Radon) integration. For simplicity, we will be interested in a particular class that we construct as follows: consider $\nu :A \longmapsto \nu(A)=\mathbb{E}(v(d\xi,A))$ where $v(d\xi,A)$ is the total variation of the random measure $d\xi$. It is defined in [\citet{Dun58} , p. 97], for all Borel set $A$ by:
\begin{equation}\label{eq2-39}
v(d\xi,A)= \sup_{I \hbox{ finite}}\left\{\sum_{i\in I}{|d\xi(A_i)|},
(A_i)_{_{i\in I}} \hbox{ partition of } A \right\}.
\end{equation}
According to \citet{Dun58}, the total variation
$v(d\xi, . )$ is a positive random measure. We deduce then that $\nu$
is a positive measure.
 We denote by $\Lambda_{\alpha}(d\xi)$ the completion, with respect to the $L^1(\nu)$ norm, of complex simple functions with bounded supports.
\section{Spectral representation of the covariation function}
\citet{Cam83} have given a spectral representation of the covariation function of  processes $X$ given in (\ref{eq2-1}) where $\xi$ have independent increments. This last condition is the key of his spectral representation. Now that we have given a weaker condition of additivity in theorem \ref {th2-1}, we generalize the Cambanis's integral representation of the covariation to processes $X$ where increments of $\xi$ are not necessarily independents but  verifying the additivity condition $\mathcal{A}$. Our spectral representation is in the sense of the general \citet{Mor56} integration with respect to the bimeasure $F$.
\subsection{Integral representation of the covariation function}
It is easy to see that for simple functions $\displaystyle f=\sum_{i=1}^{n}{f_i \1_{A_i}}$ and $\displaystyle
g=\sum_{j=1}^{m}{g_j \1_{A_j}}$ from $\Lambda_{\a}(d\xi)$:
\begin{equation}\label{eq2-40}
I_1(f,g)=I_2(f,g)=\sum_{i=1}^{n}\sum_{j=1}^{m}{f_i g_j F(A_i,A_j)}= \left[\int{f d\xi},\int{g d\xi}\right]_{\a}.
\end{equation}
Thus, simple functions of compact supports are MT-integrables. In order to generalize the property (\ref{eq2-40}) to all functions of $\Lambda_{\alpha}(d\xi)$ we will need the forthcoming proposition:
\begin{prop}\label{prop2-2}
Suppose that $d\xi$ satisfy the additivity condition $\mathcal{A}$, then we have the properties:
\begin{enumerate}
  \item For all $A\in \mathcal{B}(\mathbb{R})$ we have, $\|d\xi(A)\|_{\alpha}\leq \Psi_{\alpha}  \nu(A)$ where $\Psi_{\alpha}$ is a constant depending on $\a$. It is given by $\Psi_{\a}= \G(1-\frac{1}{\a})$. 
  \item Let $B\in \mathcal{B}(\mathbb{R})$ be a fixed Borel set. If $A$ is verifying $\nu(A)=0$ then the total variation of the complex measure $F_B$ in $A$ is null, that is $v(F_B,A)=0$. This result is also true for the measure  $\tilde{F}_A(B)$ and B fixed.
  \item Let $B\in \mathcal{B}(\mathbb{R})$ be a fixed Borel set, then for all bounded function $f\in \Lambda_{\alpha}(d\xi)$
 we have the inequality:
  \begin{equation}\label{eq2-41}
  |\int_{\mathbb{R}}{f dF_B}|\leq \Psi_{\alpha} \|d\xi(B)\|_{\alpha}^{\alpha-1} \int_{\mathbb{R}}{|f|d\nu}.
  \end{equation}
  \item Let $f\in \Lambda_{\alpha}(d\xi)$ be a fixed bounded function and denote by $G(B)=\tilde{I}(f,B)$. then $(\nu(B)=0$ implies that the total variation of $G$ vanishes i.e. $v(G,B)=0)$.
\end{enumerate}
\end{prop}
The importance of the last proposition comes from the fact that it permits to get an integral representation of the covariation function with respect to the bimeasure $F$. It transforms the covariation of two stochastic integrals with respect to $\xi$ to a double integral with respect to the bimeasure $F$ which can be seen as a covariation of increments of $\xi$.
\begin{prop}\label{prop2-3}
Let $f$ and $g$ be two bounded functions in
$\Lambda_{\alpha}(d\xi)$. The covariation of the stochastic integrals $\dy \int_{\mathbb{R}}{f d\xi}$ on $\dy
\int_{\mathbb{R}}{g d\xi}$ is given by:
\begin{equation}\label{eq2-48}
\left[\int_{\mathbb{R}}{f d\xi},\int_{\mathbb{R}}{g
d\xi}\right]_{\alpha} =
\int_{\mathbb{R}}\int_{\mathbb{R}}{f(\lambda)
\left(g(\lambda')\right)^{<\alpha-1>}F(d\lambda,d\lambda')}.
\end{equation}
The term in the right hand side is Morse Transue integral with respect to the bimeasure $F$.
\end{prop}
The equality (\ref{eq2-48}) is similar to the Cramer integral representation of the covariance of second order processes, see for instance \citet{Rao86}. This result plays a fundamental role in this work because it permits to characterize the law of the stochastic process $X$ by the bimeasure $F$ when the deterministic functions $f(t,.)$ are known. This property will be used later to classify the dependence structure of harmonisable S$\a$S processes according to the morphology of the bimeasure F.
\begin{prop}\label{prop2-4}
Suppose that  $\xi$ is a real or complex isotropic\footnote{ A complex random variable $Y$ is said isotropic if for all $\omega$, the random variables $Y$ and $Y e^{^{\imath \omega}}$ have the same distribution [\citet{Sam94}, p. 84 ]. It is known that isotropic random variables are parametric because their distribution is completely determined by their scale parameter(covariation norm).} $\a$-stable process where $d\xi$ satisfy the additivity condition $\mathcal{A}$. Then the bimeasure $F$ defined in (\ref{eq2-35}) is the unique bimeasure characterizing  $X$ and verifying the integral representation (\ref{eq2-48}).
\end{prop}

\subsection{Applications to harmonisable processes}
In this subsection we are interested in harmonisable 
S$\alpha$S processes. This means processes having the integral representation (\ref{eq2-1}) where ( $f(t,\l)=e^{\imath t\l}$ ) and  $d\xi$
 verifying the condition $\mathcal{A}$. We also suppose that the bimeasure 
$F$ given in (\ref{eq2-35}) is of bounded variation in Vitali's sense. In this case, according to \citet{Hor77} or \citet{Klu81}, $F$ can be extended to a complexe measure on the $\sigma$-algebra $\mathcal{B}(\mathbb{R}^2)$ induced by $\mathcal{B}(\mathbb{R})\times \mathcal{B}(\mathbb{R})$. By the sequel the MT-integration is the same as the usual Lebesgue (Radon) sense, see for instance \citet{Rao86}. Let us consider the covariation function denoted by $C:(s,t)\longmapsto [X_s,X_t]_{\alpha}$ that we suppose uniformly continuous, then:
\begin{itemize}
\item We will say that $X$ is \textit{\textbf{covariation stationary}} if the function $C(s,t)$ depend only on the time difference $s-t$.
\item We will say that the process  $X$ is \textit{\textbf{periodically covariated}} with a period $T>0$ if $C(s+T,t+T)=C(s,t)$ for all $s$ and $t$. 
\item We will say that the process  $X$ is \textit{\textbf{almost periodically covariated}} if $C$ is an almost periodic function in the sense of Bohr. \citet{Deh93} or \citet{Deh94}.
\end{itemize}

\begin{remark}
It should be noted that with Cambanis's spectral representation where $\xi$ have independent increments, the covariation of harmonisable processes is expressed as $\dy C(s,t)=\int_{\R}{e^{i(s-t)\l} \mu(d\l)}$ where $\mu(A)=\|d\xi(A)\|_{\a}^{\a} = [d\xi(A),d\xi(A)]_{\alpha}$ is a spectral measure on $\mathcal{B}(\R)$. We see clearly that this kind of processes are always covariation stationary. In our general spectral representation, the covariation function and hence the dependence structure of the process depend on the structure of the bimeasure $F$. Indeed, by using (\ref{eq2-48}) and the fact that ($(e^{\imath t\l'})^{<\a-1>}=e^{-\imath t\l'}$), the covariation function can be written as:
$$
C(s,t)=\int_{\R}\int_{\R}{e^{i(s\l-t\l')}F(d\l,d\l')}.
$$
The importance of our representation in the harmonisable case is that the covariation function is expressed as the Fourier transform of $F$. This result is very interesting: firstly it provides valuable informations on the link between dependence nature of the process and the morphology of its bimeasure. Secondly it makes possible statistical estimation of $F$ from the process observations. 
\end{remark}
According to \citet{Gla63}, the periodicity (resp. the almost-periodicity) of the bivariate covariation function $C$ is equivalent to the periodicity (resp. the almost-periodicity) of the univariate maps, $C_{\tau}: t\longmapsto C(t+\tau,t)$
for all $\tau$ fixed. As for periodically and almost periodically correlated processes see \citet{Hur91, Hur89}, this suggest that for all fixed $\t$ the function $C_{\tau}$ may be expressed as a Fourier series. Indeed, when $C_{\tau}$ is periodic of period $T$, 
\begin{equation}\label{eq2-51}
C_{\tau}(t) \sim \sum_{k\in \mathbb{Z}}{a_k(\tau)
e^{\imath \frac{2\pi k}{T}t}}.
\end{equation}
The convergence sense of the right hand series depends on the smoothness of $C_{\tau}$. The terms  $a_k(\t)$'s are the $k^{\hbox{th}}$ order Fourier coefficient given by, $\dy a_k(\t)=\frac{1}{T} \int_0^T{C_{\tau}(t) e^{-\imath \frac{2\pi k}{T}t} dt }$. Similarly, in the almost periodic case, the Fourier decomposition is given by:
\begin{equation}\label{eq2-51b}
C_{\tau}(t) \sim \sum_{\g\in \mathcal{J}_{\tau}}{a_{\g}(\tau)
e^{\imath \g t}},
\end{equation}
where $a_{\g}=\dy \lim_{M\rightarrow \infty}\frac{1}{M} \int_{-M}^M{C_{\tau}(t) e^{-\imath \g t} dt }$ are Bohr-Fourier  coefficients of $C_{\t}$ and $\mathcal{J}_{\tau} = \left\{\gamma \in \R, \hbox{ such that, }a_{\g}(\tau)\neq 0 \right\}$ is at most countable, see for instance \citet{Deh93}. In this last paper there was shown that according to the uniform continuity of $C$, the Bohr-Fourier coefficients $a_{\g}(\tau)$ are continuous with respect to $\tau$. This fact permit to show that $\mathcal{J} = \bigcup_{\tau}  \mathcal{J}_{\tau} $ = $\left\{\gamma \in \R, \hbox{ such that, } a_{\g}(\tau)\neq 0  \hbox{ for some } \tau\right\}$ is also countable, see \citet{Deh93}. 
\par
In the next proposition we show that some fundamental results concerning spectral representation of periodically and almost periodically correlated processes remain true for their analogues periodically and almost periodically covariated $\alpha$-stables processes. As in the works of \citet{Hur89,Hur91}, we give a necessary and sufficient condition for $X$ to be periodically or almost periodically covariated or covariation stationary.
 \begin{prop}\label{prop2-5b}
\begin{enumerate} 
\item The process $X$ is covariation stationary if and only if the bimeasure $F$ is concentrated on the diagonal line.
\item The process $X$ is periodically covariated with a period $T$ if and only if the bimeasure $F$ is concentrated on the lines parallel to the diagonal and equidistant of $\frac{2\pi}{T}$ that is the lines $(L_k)_k$ defined by,
$L_k=\{(\lambda,\lambda')\in \mathbb{R}\times \mathbb{R}, /
\lambda -\lambda'=\frac{2\pi k}{T}\}$.
\item the process $X$ is almost periodically covariated if and only if the bimeasure $F$ is supported by lines parallel to the diagonal but not necessarily equidistants i.e. the lines $(L_{\g})_{\g}$ where,
$L_{\g}=\{(\lambda,\lambda')\in \mathbb{R}\times \mathbb{R}, /
\lambda -\lambda'=\gamma \in \mathcal{J}\}$.
\end{enumerate}
\end{prop}

\bibliography{thesebibgen}

\begin{thebibliography}{20}
\providecommand{\natexlab}[1]{#1}
\providecommand{\url}[1]{\texttt{#1}}
\expandafter\ifx\csname urlstyle\endcsname\relax
  \providecommand{\doi}[1]{doi: #1}\else
  \providecommand{\doi}{doi: \begingroup \urlstyle{rm}\Url}\fi

\bibitem[Cambanis(1983)]{Cam83}
S.~Cambanis.
\newblock Complex symetric stable variables and processes.
\newblock \emph{Contribution to Statistics : P.K sen editions, Contribution to
  statistics (Essays in honor of Norman L. Johnson)}, pages 63--79, 1983.

\bibitem[Cambanis and Miamee(1989)]{Cam89}
S.~Cambanis and A.~G. Miamee.
\newblock On prediction of harmonizable stable processes.
\newblock \emph{Sankhya: The Indian Journal of Statistics, V}, 51:\penalty0
  269--294, 1989.

\bibitem[Cambanis and Miller(1980)]{Cam80}
S.~Cambanis and G.~Miller.
\newblock Some path properties of $p^{th}$ order and symmetric stable
  processes.
\newblock \emph{The Annals of Probability}, 8\penalty0 (6):\penalty0
  1148--1156, 1980.

\bibitem[Cambanis and Miller(1981)]{Cam81}
S.~Cambanis and G.~Miller.
\newblock Linear problems in $p^{th}$ order and symmetric stable processes.
\newblock \emph{SIAM Journal of Applied Mathematics}, 41:\penalty0 43--69,
  1981.

\bibitem[Cambanis and Soltani(1984)]{Cam84}
S.~Cambanis and A.~R. Soltani.
\newblock Prediction of stable processes: Spectral and moving average
  representation.
\newblock \emph{Zeitschrift fäur Wahrscheinlichkeitstheorie und verwandte
  Gebiete}, 66:\penalty0 593--612, 1984.

\bibitem[Chang and Rao(1986)]{Rao86}
DK~Chang and MM~Rao.
\newblock {Bimeasures and nonstationary processes, in " Real and stochastic
  analysis"}.
\newblock John Wiley and Sons, New York, 1986.

\bibitem[Dehay and Hurd(1993)]{Deh93}
D.~Dehay and H.~L. Hurd.
\newblock Representation and estimation for periodically and almost
  periodically correlated random processes.
\newblock \emph{Cyclostationarity in Communications and Signal Processing},
  1993.

\bibitem[Dehay(1994)]{Deh94}
Dominique Dehay.
\newblock Spectral analysis of the covariance of the almost periodically
  correlated processes.
\newblock \emph{Stochastic processes and their applications}, 50\penalty0
  (2):\penalty0 315--330, 1994.

\bibitem[Dunford and Schwartz(1958)]{Dun58}
N.~Dunford and J.~T. Schwartz.
\newblock \emph{Linear Operators, Part I}.
\newblock Wiley, New York, 1958.

\bibitem[Gladyshev(1963)]{Gla63}
E.~G. Gladyshev.
\newblock Periodically and almost periodically correlated random processes with
  continuous time parameter.
\newblock \emph{Theory Prob. and Appl}, 8:\penalty0 173--177, 1963.

\bibitem[Horowitz(1977)]{Hor77}
J.~Horowitz.
\newblock Une remarque sur les bimesures.
\newblock \emph{S{\'e}minaire de probabilit{\'e}s de Strasbourg}, 11:\penalty0
  59--64, 1977.

\bibitem[Hurd(1989)]{Hur89}
H.L. Hurd.
\newblock Representation of strongly harmonizable periodically correlated
  processes and their covariances.
\newblock \emph{Journal of Multivariate Analysis}, 29\penalty0 (1):\penalty0
  53--67, 1989.

\bibitem[Hurd(1991)]{Hur91}
H.L. Hurd.
\newblock Correlation theory of almost periodically correlated processes.
\newblock \emph{Journal of Multivariate Analysis}, 37\penalty0 (1):\penalty0
  24--45, 1991.

\bibitem[Kluvanek(1981)]{Klu81}
I.~Kluvanek.
\newblock Remarks on bimeasures.
\newblock \emph{Proceedings of the American Mathematical Society}, 81\penalty0
  (2):\penalty0 233--239, 1981.

\bibitem[Miller(1977)]{Mil77}
G.~W. Miller.
\newblock \emph{Some results on symmetric stable distributions and processes}.
\newblock PhD thesis, University of North Carolina at Chapel Hill., 1977.

\bibitem[Morse and Transue(1956)]{Mor56}
M.~Morse and W.~Transue.
\newblock C-bimeasures $\lambda$ and their integral extensions.
\newblock \emph{The Annals of Mathematics}, 64\penalty0 (3):\penalty0 480--504,
  1956.

\bibitem[Nikias and Shao(1996)]{Nik96}
C.~L. Nikias and M.~Shao.
\newblock \emph{Signal processing with alpha-stable distributions and
  applications}.
\newblock New York :Wiley, 1996.

\bibitem[Rao(1982)]{Rao82}
M.~M. Rao.
\newblock Harmonizable processes : Structure theory.
\newblock \emph{L'enseignement Mathématiques (Essays in Honor of prof. S.
  Bochner)}, 28:\penalty0 295--351, 1982.

\bibitem[Samorodnitsky and Taqqu(1994)]{Sam94}
G.~Samorodnitsky and M.~S. Taqqu.
\newblock \emph{Stable non-Gaussian random processes: stochastic models with
  infinite variance}.
\newblock Chapman \& Hall, 1994.

\bibitem[Uchaikin and Zolotarev(1999)]{Uch99}
V.~V. Uchaikin and V.~M. Zolotarev.
\newblock \emph{Chance and Stability: Stable Distributions and Their
  Applications}.
\newblock VSP, 1999.

\end{thebibliography}
\newpage
\section{Appendix}

{\bf Proof of theorem \ref{th2-1}: }

We begin with the case when $X^d$ is an 
S$\a$S real vector and then we generalize to the complex case.
\\
\underline{\textbf{Case when $X^d$ is real:}}

We begin with demonstrating the additivity in the case of three components. By definition of the covariation and by  using the result (\ref{eq2-3}) it is easy to see that, for all $i_0, j$
and $k$ such that $j\neq k$ and for all reals $\theta_j$ and
$\theta_k$ we have:
\begin{equation}\label{eq2-21}
\begin{array}{lll}
D & \triangleq & \displaystyle
[X_{i_0},\th_j X_j+\th_k X_k]_{_{\alpha}}-[X_{i_0},\th_j X_j]_{_{\alpha}}-[X_{i_0},\th_k X_k]_{_{\alpha}},\\
& = &\displaystyle \int_{\S_{d}}{s_{i_0}  \left((\th_j s_j+
\th_k.s_k)^{^{<\alpha-1>}}-(\th_j
s_j)^{^{<\alpha-1>}}-(\th_k s_k)^{^{<\alpha-1>}}\right)d\Gamma_{X^{d}}},\\
& = & \displaystyle \int_{\S_{d}}{s_{i_0}
\Delta(s_j,s_k)d\Gamma_{X^{d}}},
 \end{array}
\end{equation}
where $\Delta(s_j,s_k)=(\th_j s_j+
\th_k s_k)^{^{<\alpha-1>}}-(\th_j
s_j)^{^{<\alpha-1>}}-(\th_k s_k)^{^{<\alpha-1>}}$. For the rest of the proof we will need the next lemma:
\begin{lem} \label{lem2-1}
Let $1<p<2$, then for all real $s$, we have:
\begin{equation}\label{eq2-4}
s^{<p-1>}=\frac{1}{\Gamma(1-p) \cos(\frac{p\pi}{2})}
\int_{0}^{+\infty}{\frac{\sin(st)}{t^{p}}}dt\triangleq\varrho_{p}
\int_{0}^{+\infty}{\frac{\sin(st)}{t^{p}}}dt,
\end{equation}
where $\Gamma$ is the usual gamma function.
\end{lem}
 With the equality (\ref{eq2-4}) we deduce:
\begin{equation}\label{eq2-22}
\Delta(s_j,s_k)=\varrho_{\a}
\int_{0}^{+\infty}{\frac{\sin((\th_j s_j+\th_k s_k)t)-\sin(\th_j s_j t)-\sin(\th_k s_k t)}{t^{\a}}}dt.
\end{equation}
By using the classical trigonometric properties,
\begin{equation}\label{eq2-23}
\begin{array}{ccc}
\sin(p)-\sin(q) = 2 \cos(\frac{p+q}{2})  \sin(\frac{p-q}{2}),\\
\cos(p)-\cos(q) = -2 \sin(\frac{p+q}{2})  \sin(\frac{p-q}{2}),
\end{array}
\end{equation}
it is easy to see that:
\begin{equation}\label{eq2-24}
\Delta(s_j,s_k)=-4\varrho_{\a}
\int_{0}^{+\infty}{\frac{\sin(\frac{\th_j s_j}{2}t) \sin(\frac{\th_ks_k}{2}t) \sin(\frac{\th_js_j+\th_ks_k}{2}t)}{t^{\a}}}dt.
\end{equation}
We replace $\Delta(s_j,s_k)$ in (\ref{eq2-21}) and apply the Fubini theorem, which is applicable in our case, we have:
\begin{equation}\label{eq2-24b}
\begin{array}{llll}
 D  & = &- \displaystyle 4\varrho_{\a}\int_{0}^{+\infty}{\frac{\dy \int_{\S_{d}}{s_{i_0}
\sin(\frac{\th_js_j}{2}t) \sin(\frac{\th_ks_k}{2}t) \sin(\frac{\th_js_j+\th_ks_k}{2}t)d\Gamma_{X^{d}}}}{t^{\a}}dt},\\
& = & - \displaystyle 4\varrho_{\a} \int_{0}^{+\infty}{\frac{S(t
\th_j,t \th_k)}{t^{\a}}dt},
\end{array}
\end{equation}
where $S(x,y)=\dy \int_{\S_{d}}{s_{i_0}
\sin(\frac{xs_j}{2}) \sin(\frac{ys_k}{2}) \sin(\frac{xs_j+ys_k}{2})d\Gamma_{X^{d}}}$.
Since the unit sphere $\S_d$ is compact, then the measure $\Gamma_{_{X^d}}$ have finite moments of all orders. By the sequel the  function $S(x,y)$ is infinitely differentiable. Its partial derivative withe respect to $y$ is given by:
\begin{equation}\label{eq2-25}
\frac{\partial S}{\partial y}(x,y)=\frac{1}{2}\int_{\S_{d}}{s_{i_0} s_k
\sin(\frac{xs_j}{2}) \sin(\frac{xs_j+2ys_k}{2})d\Gamma_{X^{d}}}.
\end{equation}
By differentiating the above equality, with respect to  $x$, we have:
\begin{equation}\label{eq2-26}
\frac{\partial^2 S}{\partial x
\partial y}(x,y)=\frac{1}{4}\int_{\S_{d}}{s_{i_0} s_j s_k
 \sin((x s_j + y s_k))d\Gamma_{X^{d}}}.
\end{equation}
Let us remark by using (\ref{eq2-18}) that:
$$
\frac{\partial^2 S}{\partial x
\partial
y}(\th_j,\th_k)=\frac{1}{4}\frac{\partial^3\phi}{\partial\theta_{i_0}\partial\theta_j\partial\theta_k}(0,\dots,0,\th_j,0,\dots,0,\th_k,0,\dots,0).
$$
Therefore, according to the condition (\ref{eq2-20}),
we have:
\begin{equation}\label{eq2-27}
\forall \th_j, \th_k \in \mathbb{R}\hspace{1cm} \frac{\partial^2
S}{\partial x
\partial y}(\th_j,\th_k)=0.
\end{equation}
By (\ref{eq2-24b}) and (\ref{eq2-27}), in order to show that,
$\displaystyle
[X_i,\th_j X_j+\th_k X_k]_{_{\alpha}}=[X_i,\th_j X_j]_{_{\alpha}}+[X_i,\th_k X_k]_{_{\alpha}}$,
it is sufficient to show that (\ref{eq2-27}) is equivalent to the fact that
$S$ is a vanishing function, that is:
\begin{equation}
\forall x,y, \frac{\partial^2 S}{\partial x\partial y}(x,y)=0
\Longleftrightarrow \forall x,y,  S(x,y)=0.
\end{equation}
The second implication is obvious because if $S$ is a null function, it is the same for its second derivative. Now assume that we have, $\dy  \forall x,y, \frac{\partial^2
S}{\partial x\partial y}(x,y)=0$. We deduce then that $\dy \forall
x,y, \frac{\partial S}{\partial y}(x,y)$ does not depend on $x$.
This implies that, $\dy \forall x,y, \frac{\partial
S}{\partial y}(x,y)=\frac{\partial S}{\partial y}(0,y)$. By replacing the right term of this last equality in the formula(\ref{eq2-25}), we find that $\dy \forall x,y,
\frac{\partial S}{\partial y}(x,y)=0$. In the same way, this implies that
$S(x,y)$ does not depend on $y$, consequently $\ \forall x,y,
S(x,y)=S(x,0)$ and by replacing this in the formula of $S$, given in (\ref{eq2-24b}), we see that $S$ is null function.
\par
The proof of the general case is done with a technique
similar to the one of three variables discussed above. Indeed,
let $\th_1, \th_2,\dots,\th_d$ be real numbers and fix
$i_0\in\{1,..,d\}$. We note $D$ the term expressed as:
\begin{equation}\label{eq2-28}
\begin{array}{lll}
D& \triangleq & \displaystyle
[X_{_{i_0}},\th_{_1}X_{_1}+\dots+\th_{_d}X_{_d}]_{_{\alpha}}-[X_{_{i_0}},\th_{_1}X_{_1}]_{_{\alpha}}-\dots-[X_{_{i_0}},\th_{_d}
X_{_d}]_{_{\alpha}},\\
& = & \displaystyle
[X_{_{i_0}},\th_{_1}X_{_1}+\dots+\th_{_d}X_{_d}]_{_{\alpha}}-[X_{_{i_0}},\th_{_1}X_{_1}+\dots+\th_{_{d-1}}X_{_{d-1}}]_{_{\alpha}}-[X_{_{i_0}},\th_{_d}
X_{_d}]_{_{\alpha}},\\
& + & \displaystyle
[X_{_{i_0}},\th_{_1}X_{_1}+\dots+\th_{_{d-1}}X_{_{d-1}}]_{_{\alpha}}-[X_{_{i_0}},\th_{_1}X_{_1}+\dots+\th_{_{d-2}}X_{_{d-2}}]_{_{\alpha}}-[X_{_{i_0}},\th_{_{d-1}}
X_{_d}]_{_{\alpha}},\\
& + & \displaystyle \dots+
[X_{_{i_0}},\th_{_1}X_{_1}+\th_{_2}X_{_{2}}]_{_{\alpha}}-[X_{_{i_0}},\th_{_1}X_{_1}]_{_{\alpha}}-[X_{_{i_0}},\th_{_{2}}
X_{_2}]_{_{\alpha}},\\
& = & \displaystyle \sum_{j=0}^{d-2}{D_j}.
\end{array}
\end{equation}
For all  $j\in\{0,\dots,d-2\}$, the terms $D_j$'s are given by:
\begin{equation}\label{eq2-29}
\begin{array}{lll}
D_j & \triangleq
&[X_{_{i_0}},\th_{_1}X_{_1}+\dots+\th_{_{d-j}}X_{_{d-j}}]_{_{\alpha}}-[X_{_{i_0}},\th_{_1}X_{_1}+\dots+\th_{_{d-j-1}}X_{_{d-j-1}}]_{_{\alpha}}-[X_{i_0},\th_{d-j}X_{d-j}
]_{_{\alpha}},\\
& = & \displaystyle \int_{\S_{d}}s_{i_0}\left[( \th _1 s_1+\dots+\th_
{_{d-j}} s_{_{d-j}})^{^{<\alpha-1>}}\right.\\
&  & \hspace{2.5cm} \dy\left. -( \th _1 s_1+\dots+\th_{_{d-j-1}}
s_{_{d-j-1}})^{^{<\alpha-1>}}
 -(\th_{_{d-j}}s_{_{d-j}})^{^{<\alpha-1>}}\right]d\Gamma_{X^{d}}.
 \end{array}
\end{equation}
Like in the case of three variables discussed above, we use the equality (\ref{eq2-4}) of lemma \ref{lem2-1} and then we apply the Fubini theorem. We find an equality similar to  (\ref{eq2-21}). It is given by:
\begin{equation}\label{eq2-3O}
\begin{array}{lll}
D_j & = & \displaystyle
\varrho_{\a}\int_{0}^{+\infty}{\frac{\int_{\S_{d}}{s_{i_0}
\sin(\frac{\th_{d-j}s_{d-j}}{2}t) \sin(\frac{\th_1
s_1+\dots+\th_{_{d-j-1}} s_{_{d-j-1}}}{2}t).\sin(\frac{\th_1
s_1+\dots+\th_{_{d-j}} s_{_{d-j}}}{2}t)d\Gamma_{X^{d}}}}{t^{\a}}dt},\\
& = & \displaystyle
\varrho_{\a}\int_{0}^{+\infty}{\frac{\Delta_j(t \th_1,\dots,t
\th_{d-j})}{t^{\a}}dt},
 \end{array}
\end{equation}
where,
\begin{equation}\label{eq2-31}
\begin{array}{lll}
\dy \Delta_j(\th_{_1},\dots,\th_{_{d-j}})& \triangleq &\dy
\int_{S_{_{d}}}{s_{_{i_0}}
\sin\left(\frac{\th_{_{d-j}}s_{_{d-j}}}{2}\right) \sin\left(\frac{\th
_{_1} s_{_1}+\dots+\th_{_{d-j-1}}
s_{_{d-j-1}}}{2}\right)}\\
&  & \hspace{3.5cm} \dy  {\sin\left(\frac{\th_{_1}
s_{_1}+\dots+\th_{_{d-j}} s_{_{d-j}}}{2}\right)d\Gamma_{_{X^{d}}}}.
\end{array}
\end{equation}
First, the functions
$\Delta_j$'s are two-differentiables because the measure $\Gamma_{X^{d}}$ have finite moment of all orders. With the same technique as in (\ref{eq2-25}) and (\ref{eq2-26}), for all
$k=1,\dots,d-j-1$, by differentiating the function  $\Delta_j$ with respect to
$\th_k$ and $\th_{d-j}$, we have:
\begin{equation}\label{eq2-32}
\begin{array}{lll}
\displaystyle  \frac{\partial^2 \Delta_j}{\partial \th_k
\partial \th_{d-j}}(\th_1,\dots,\th_{d-j})& = & \displaystyle  \frac{1}{4}\int_{\S_{d}}{s_{i_0} s_k s_{d-j}
 \sin(\frac{\th_1 s_1+\dots2 \th_k s_k+\dots+2 \th_{d-j} s_{d-j}}{2})d\Gamma_{X^{d}}},\\
& = &\displaystyle
\frac{1}{4}\frac{\partial^3\phi}{\partial\th_{i_0}\partial \th_k
\partial\theta_{d-j}}(\frac{\th_1}{2},\dots, \th_k,\frac{\th_{k+1}}{2},\dots,\th_{d-j},0,\dots,0).
\end{array}
\end{equation}
By using  (\ref{eq2-18}) and the condition (\ref{eq2-20}), we deduce that:
\begin{equation}\label{eq2-33}
 \forall k=1,\dots,d-j-1,
\forall \th_1,\dots,\th_{d-j},  \frac{\partial^2 \Delta_j}{\partial
\th_k \partial \th_{d-j}}(\th_1,\dots,\th_{d-j})=0.
\end{equation}
According to (\ref{eq2-3O}) and (\ref{eq2-33}), in order to show that $D=0$, 
it suffices to show the next equivalence:
$$
\begin{array}{ccc}
\dy \forall k=1,\dots,d-j-1,  \forall \th_1,\dots,\th_{d-j},
 \frac{\partial^2, \dy
\Delta_j}{\partial \th_k\partial
\th_{d-j}}(\th_1,\dots,\th_{d-j})=0 \\\Longleftrightarrow\\
\dy \forall \th_1,\dots,\th_{d-j}, \Delta_j(\th_1,\dots,\th_{d-j})=0.
\end{array}
$$
A trivial implication is that if $\Delta_j$ is a null function, it is the same for its second derivative. For the other implication, let us assume that $\forall k=1,\dots,d-j-1$ and $\forall
\th_1,\dots,\th_{d-j}$ we have, $\dy  \frac{\partial^2
\Delta_j}{\partial \th_k\partial
\th_{d-j}}(\th_1,\dots,\th_{d-j})=0$. One concludes, therefore, that $\forall
\th_1,\dots,\th_{d-j}$, the function $\dy \frac{\partial
\Delta_j}{\partial \th_{d-j}}(\th_1,\dots,\th_{d-j})$ does not depend on $\th_1,\dots,\th_{d-j-1}$, which allows us to write, $ \forall
\th_1,\dots,\th_{d-j}$, $\dy \frac{\partial \Delta_j}{\partial
\th_{d-j}}(\th_1,\dots,\th_{d-j})=\frac{\partial \Delta_j}{\partial
\th_{d-j}}(0,\dots,0,\th_{d-j})$. By replacing the right term of this last equality in the formula  $\dy \frac{\partial
\Delta_j}{\partial \th_{d-j}}$, we find that $\forall \th_1,\dots,\th_{d-j}$, $\dy \frac{\partial
\Delta_j}{\partial \theta_{d-j}}(\th_1,\dots,\th_{d-j})=0$. It also implies that $\Delta_j$ does not depend on $\th_{d-j}$, consequently $\forall$, $ \th_1,\dots,\th_{d-j}$, we have the equality 
$\Delta_j(\th_1,\dots,\th_{d-j})=\Delta_j(\th_1,\dots,\th_{d-j-1},0)$.
It is sufficient to replace $\theta_{d-j}$ by 0 in the formula (\ref{eq2-31}), to see that  $\Delta_j$ is a vanishing function.

\underline{\textbf{Case where $X^d$ is a complex vector:}}

\par
The prof of theorem \ref{th2-1} when $X^d$ is complex is similar to the real case. We give an idea on how to show this result in the case of three variables. First by definition of the covariation, for all $i_{0}, j$ and $k$ not all equals and for all complex numbers $\theta_j$ and $\theta_k$ we have:
\begin{equation}\label{eq2-21bb}
\begin{array}{lll}
D & \triangleq & \displaystyle
[X_{i_0},\th_j X_j+\th_k X_k]_{_{\alpha}}-[X_{i_0},\th_j X_j]_{_{\alpha}}-[X_{i_0},\th_k X_k]_{_{\alpha}},\\
& = &\displaystyle \int_{\S_{2d}}{s_{i_0}  \left((\th_j s_j+
\th_k s_k)^{^{<\alpha-1>}}-(\th_j
s_j)^{^{<\alpha-1>}}-(\th_k s_k)^{^{<\alpha-1>}}\right)d\Gamma_{X^{2d}}},\\
& = & \displaystyle \int_{\S_{2d}}{s_{i_0}
\Delta(s_j,s_k)d\Gamma_{X^{2d}}},
 \end{array}
 \end{equation}
where $\Delta(s_j,s_k)=(\th_j s_j+ \th_k s_k)^{^{<\alpha-1>}}-(\th_j s_j)^{^{<\alpha-1>}}-(\th_k s_k)^{^{<\alpha-1>}}$. Instead of using lemma \ref{lem2-1}, we use the following result:
 \begin{lem}\label{lem2-2}
Let $z=a+\imath b$ be complex number and $p$ positifs number such that $0<p<2$. Then by using the complex notation $x=s+\imath t$, we have the following equality:
\begin{equation}\label{eq2-6}
\int_{-\infty}^{\infty}{\int_{-\infty}^{\infty}{\frac{1-\cos(\RE(x \overline{z}))}{|x|^{(p+2)}}ds}dt}=c(p)
 |z|^{p},
\end{equation}
 where $\dy c(p)=2^{\frac{-p}{2}}\frac{\Gamma(1-p) \cos(\frac{p\pi}{2})}{p} \int_{0}^{2
\pi}{\left|1+\sin(2\theta)\right|^{\frac{p}{2}}d\theta}
$. We have a result analogous to (\ref{eq2-4}):
\begin{equation}\label{eq2-9}
 z^{<p-1>}= p c(p) \int_{-\infty}^{\infty}{\int_{-\infty}^{\infty}{\frac{\sin(\RE(x \overline{z}))}{(\overline{x})^{<p+1>}}ds}dt}\triangleq \rho_p .
\int_{-\infty}^{\infty}{\int_{-\infty}^{\infty}{\frac{\sin(\RE(x \overline{z}))}{(\overline{x})^{<p+1>}}ds}dt}.
\end{equation}
 \end{lem}
In this case, by applying the equality (\ref{eq2-9}), for all $x=x_1+
\imath x_2$ we have :
\begin{equation}\label{eq2-22bb}
\Delta(s_i,s_j)=\rho_{\a} \int\int_{-\infty}^{+\infty}{ \frac{\sin(\RE(\th_j s_j+\th_k s_k)\overline{x})-\sin(\RE(\th_js_j\overline{x}))-\sin(\RE(\th_ks_k\overline{x}))}{(x)^{<\a>}}dx_1
dx_2}.
\end{equation}
The rest of the proof is the same as in the case of real S$\a$S vectors, where  we use the lemma \ref{lem2-2}, the trigonometric properties (\ref{eq2-23}) and the complex derivative operator as in (\ref{eq2-32}) and (\ref{eq2-33}). Finally the same reasoning as in the real case shows the results.

{\bf Proof of proposition \ref{prop2-2} : }

\begin{enumerate}
\item Let $A\in \mathcal{B}( \mathbb{R} )$ be  a fixed Borel set. Since $d\xi(A)$ is a real (resp. isotropic complex) S$\a$S random variable then according to \citet{Cam89} or \citet{Sam94}, the $p^{\hbox{th}}$-order fractional moments of $d\xi(A)$ are given by $\EE|d\xi(A)|^p=S_{\a}(p) \|d\xi(A)\|_{\alpha}^{p}$ where $S_{\a}(p) = 2^{p}
\frac{\G(\frac{1+p}{2}) \G(1-\frac{p}{\a}) }{\G(1-\frac{p}{2}) \G(\frac{1}{2})} $ (resp. $S_{\a}(p) = \frac{\G(\frac{2+p}{2}) \G(1-\frac{p}{\a}) }{\G(1-\frac{p}{2})}$). In particular, for $p=1$, $S_{\a}(1)\leq \Psi_{\a}$ then  $\|d\xi(A)\|_{\alpha} \leq \Psi_{\a}.  \mathbb{E}|d\xi(A)| $. According to the definition of the total variation, defined in (\ref{eq2-39}), it is easy to see that $|d\xi(A)| \leq
v(d\xi,A)$. Therefore, we have:
\begin{equation}\label{eq2-42}
\|d\xi(A)\|_{\alpha}=\Psi_{\a} \mathbb{E}|d\xi(A)|\leq
\Psi_{\a} \mathbb{E}\left(v(d\xi,A)\right)=
\Psi_{\a} \nu(A).
\end{equation}
 \item Let $B\in \mathcal{B}(\mathbb{R})$ be fixed Borel set, then by definition of the total variation of the complex measure $F_B$ given in 
(\ref{eq2-39}) (see \citet{Dun58}) and replacing 
$F_B$ by its expression we have:
\begin{displaymath}
\begin{array}{llll}
\displaystyle v(F_B,A)& = & \displaystyle \sup_{I\hbox{ finite}}\left\{\sum_{i\in I}{|F_B(A_i)|}, (A_i)_{i\in I} \hbox{ partition of A} \right\}, \\
& = & \displaystyle \sup_{I\hbox{ finite}}\left\{\sum_{i\in
I}{|[d\xi(A_i),d\xi(B)]_{\alpha}|}, (A_i)_{i\in I} \hbox{
partition of A} \right\}.
\end{array}
\end{displaymath}
It is known that $|[d\xi(A_i),d\xi(B)]_{\alpha}|\leq \|d\xi(A_i)\|_{\alpha} \|d\xi(A_i)\|^{\a-1}_{\alpha}$, see for instance [\citet{Sam94}, p.96]. Then by using (\ref{eq2-42}) we deduce that,
\begin{displaymath}
\begin{array}{llll}
\displaystyle v(F_B,A) & \leq & \displaystyle
\|d\xi(B)\|_{\alpha}^{\alpha-1} \sup_{I\hbox{
finite}}\left\{\sum_{i\in
I}{\|d\xi(A_i)\|_{\alpha}}, (A_i)_{i\in I} \hbox{ partition of A} \right\},\\
& \leq & \displaystyle \Psi_{\a}
\|d\xi(B)\|_{\alpha}^{\alpha-1}\sup_{I\hbox{
finite}}\left\{\sum_{i\in
I}{\nu(A_i)}, (A_i)_{i\in I} \hbox{ partition of A} \right\}, \\
& \leq & \displaystyle \Psi_{\a}
\|d\xi(B)\|_{\alpha}^{\alpha-1}\nu(A).
\end{array}
\end{displaymath}
By the sequel, if $\nu(A)$ is null then $v(F_B,A)=0$.
\item First we begin by showing the inequality 
(\ref{eq2-41}) for simple functions. Indeed, let $f$ be a complex simple function defined as
$\dy f=\sum_{i=1}^{n}{f_i \1_{A_i}}$, then :
\begin{equation}\label{eq2-43}
\begin{array}{lll}
\dy\left|\int_{\mathbb{R}}{f d F_B}\right|& = & \dy \left|\sum_{i}^n{f_i F_B(A_i)}\right|=\left|\sum^n _i{f_i \left[d\xi(A_i),d\xi(B)\right]_{\alpha}}\right|=\left|\left[\sum^n _i{f_i d\xi(A_i)},d\xi(B)\right]_{\alpha}\right|,\\
& \leq & \dy \left\|d\xi(B)\right\|_{\alpha}^{\alpha-1} \left\|\sum_i^n{f_i d\xi(A_i)}\right\|_{\alpha}, \\
& = & \dy \Psi_{\a} \left\|d\xi(B)\right\|_{\alpha}^{\alpha-1} .  \mathbb{E}\left|\sum_i^n{f_i d\xi(A_i)}\right|, \\
& \leq & \dy \Psi_{\a} \left\|d\xi(B)\right\|_{\alpha}^{\alpha-1} \sum_i^n{|f_i| . \mathbb{E}\left|d\xi(A_i)\right|},\\
& \leq & \dy \Psi_{\a} \left\|d\xi(B)\right\|_{\alpha}^{\alpha-1} \sum_i^n{|f_i| . \mathbb{E}\left(v(d\xi,A_i)\right)} , \hbox{ using (\ref{eq2-42})}\\
& = & \displaystyle \Psi_{\a}
\left\|d\xi(B)\right\|_{\alpha}^{\alpha-1}
\int_{\mathbb{R}}{|f|d\nu}.
\end{array}
\end{equation}
If now $f$ is a bounded function in $\Lambda_{\alpha}(\xi)$, then there exist as  sequence of  simple functions $(f_n)$ such that $\int{|f_n-f|d\nu}$ converge to 0 as $n$ tends toward infinity. First let us remark that the sequence
$(\int_{\mathbb{R}}{f_n}dF_B)_n$ is convergent because it verifies the Cauchy condition. Indeed, according to the inequality (\ref{eq2-43}), for all integers $m$ and $n$,
\begin{equation}\label{eq2-44}
\left|\int_{\mathbb{R}}{f_n}dF_B-\int_{\mathbb{R}}{f_m}dF_B\right|
\leq \Psi_{\a} \left\|d\xi(B)\right\|_{\alpha}^{\alpha-1}
\int_{\mathbb{R}}{|f_n-f_m|d\nu}.
\end{equation}
On the other hand, since $(f_n)_n$ converges to $f$ in $L_1(\nu)$
we can then extract a subsequence $(f_{n_k})_k$ that converges to $f$ ($\nu$-almost everywhere). This means that there exists a Borel set $D\in B({\mathbb{R}})$ such that $\nu(\mathbb{R}\setminus
D)=0$ and $f_{n_k}\1_{D}$ converges toward $f.\1_{D}$ as
$k$ tends to infinity. Using the second result of these proposition and the fact that $\nu(\mathbb{R}\setminus D)=0$ we deduce that $v(F_B , \mathbb{R}\setminus D)=0$, this implies that the subsequence $(f_{n_k})_k$ converges $v(F_B,.)$-almost everywhere to $f$. The functions $f_{n_k}$ and the measure $v(F_B,.)$
are bounded, then by applying the Lebesgue convergence theorem in the case of Dunford-Schwartz integrals [\citet{Dun58},  p.151], we deduce then that the subsequence $(\int_{\mathbb{R}}{f_{n_k}dF_B})_k$ converge to
$\int_{\mathbb{R}}{f dF_B}$. Using the uniqueness of the limit we deduce that the sequence $(\int_{\mathbb{R}}{f_n}dF_B)_n$ converge
to $\int_{\mathbb{R}}{f dF_B}$. Finally, by using
(\ref{eq2-43}) we have:
\begin{equation}\label{eq2-45}
|\int_{\mathbb{R}}{f_n}dF_B | \leq \Psi_{\a}
\left\|d\xi(B)\right\|_{\alpha}^{\alpha-1}
\int_{\mathbb{R}}{|f_n|d\nu}.
\end{equation}
The inequality (\ref{eq2-41}) is then obtained by  passage to the limit in   (\ref{eq2-45}).
 \item Now let $f\in\Lambda_{\alpha}(d\xi)$ be a fixed function. Using the definition of the total variation of $G$,
 for all $B\in \mathcal{B}(\mathbb{R})$, we have :
\begin{equation}\label{eq2-46}
\begin{array}{llll}
v(G,B) & = &\dy \sup_{I \hbox{ finite}}\left\{\sum_{i\in I}{|G(B_i)|}, (B_i)_{i\in I} \hbox{ partition of B} \right\},\\
& = & \displaystyle \sup_{I \hbox{ finite}}\left\{\sum_{i\in I}{|\int_{\mathbb{R}}{f dF_{B_i}|}}, (B_i)_{i\in I} \hbox{ partition of B} \right\},\\
 & \leq & \displaystyle \Psi_{\a}  \int_{\mathbb{R}}{|f| d\nu}  \sup_{I \hbox{ finite}}\left\{\sum_{i\in I}{\|d\xi(B_i)\|_{\alpha}^{\alpha-1}}, (B_i)_{i\in I} \hbox{ partition of B}
 \right\}.
\end{array}
\end{equation}
By using the inequality (\ref{eq2-42}) in (\ref{eq2-46}) we deduce that:
\begin{equation}\label{eq2-47}
v(G,B) \leq  (\Psi_{\a})^{\a} 
\int_{\mathbb{R}}{|f| d\nu} . \sup_{I \hbox{
finite}}\left\{\sum_{i\in I}{(\nu(B_i))^{\alpha-1}}, (B_i)_{i\in I}
\hbox{ partition of B} \right\}.
\end{equation}
If $\nu(B)=0$ then for all finite index set $I$ and for all partition $(B_i)_{i\in I}$ of $B$, the $B_i$'s are subsets of $B$, therefore $\nu(B_i)=0$  for all $i$ and we have $\displaystyle \sup_{I \hbox{ finite }} \left\{\sum_{i\in I}{|\nu(B_i)|^{\alpha-1}}, (B_i)_{i\in I} \hbox{ partition of B}
\right\}=0$ which implies that  $v(G,B)=0$.
\end{enumerate}

{\bf Proof of proposition \ref{prop2-3}: }

Let $f$ and $g$ be two bounded functions of $\dy \Lambda_{\alpha}(d\xi)$ and $\dy
f_n=\sum_{i}{f^{n}_{i} \1_{A^{n}_{i}}}$ and $\dy
g_n=\sum_{j}{g^{n}_{j} \1_{B^{n}_{j}}}$ two sequences of simple functions in  $\Lambda_{\alpha}(d\xi)$ such that
$\int_{\mathbb{R}}{|f_n-f|d\nu}$ and
$\int_{\mathbb{R}}{|g_n-g|d\nu}$
converge to 0. The equality (\ref{eq2-48}) is satisfied for the sequences of simple functions $(f_n)_n$ and $(g_n)_n$. Indeed, let $m$ and $n$ be two integers, then:
$$
\begin{array}{lll}
\dy \left[\int_{\mathbb{R}}{f_n d\xi},\int_{\mathbb{R}}{g_m
d\xi}\right]_{\alpha} & = & \dy \left[\sum_{i}{f^{n}_{i}
d\xi(A^{n}_{i})},\sum_{j}{g^{n}_{j}
d\xi(B^{n}_{j})}\right]_{\alpha},\\
 & = & \dy \sum_{i}{f^{n}_{i} \left[
d\xi(A^{n}_{i}),\sum_{j}{g^{n}_{j}
d\xi(B^{n}_{j})}\right]_{\alpha}}.
\end{array}
$$
The Borel sets $(B_j^n)_j$ are pairwise disjoints and since
$d\xi$ satisfy the condition ($\mathcal{A}$) we deduce that:
\begin{equation}\label{eq2-49}
\begin{array}{lll}
\dy \left[\int_{\mathbb{R}}{f_n d\xi},\int_{\mathbb{R}}{g_m
d\xi}\right]_{\alpha}& = & \dy  \sum_{i}\sum_{j}{f^{n}_i \left(g^{m}_j\right)^{<\alpha-1>} [d\xi(A^{n}_i),d\xi(B^{m}_j)]_{\alpha}},\\
&=& \dy \sum_{i}\sum_{j}{f^{n}_i \left(g^{m}_j\right)^{<\alpha-1>} F(A^{n}_i,B^{m}_j)},\\
&=& \dy \int_{\mathbb{R}}\int_{\mathbb{R}}{f_{n}(\lambda)
\left(g_{m}(\lambda')\right)^{<\alpha-1>}F(d\lambda,d\lambda')}.
\end{array}
\end{equation}
The result (\ref{eq2-49}) remain true for any function $f$ and $g$ in $\Lambda_{\alpha}(d\xi)$. First, remark that:
\begin{displaymath}
\sum_{i}\sum_{j}{f^{n}_i \left(g^{m}_j\right)^{<\alpha-1>} [d\xi(A^{n}_i),d\xi(B^{m}_j)]_{\alpha}}\\
=\sum_{j}\left(g^{m}_j\right)^{<\alpha-1>}
\int{f_n(\lambda)F_{B^{m}_j}(d\lambda)}.
\end{displaymath}
Firstly, for a fixed $m$ and using the same reasoning as in the proof of the third assertion of the proposition
\ref{prop2-2}, we have: $\dy
\int{f_n(\lambda)F_{B^{m}_j}(d\lambda)}$ converge to $\dy
\int{f(\lambda)F_{B^{m}_j}(d\lambda)}$. On the other hand, using (\ref{eq2-42}) we can derive $\left\|\int_{\mathbb{R}}{f_n-f d\xi}\right\|_{\a}\leq \Psi_{\a} \int{|f_n-f|d\nu}$. We have then the inequality,
$$
\begin{array}{lll}
\dy\left| \left[\int_{\mathbb{R}}{f_n d\xi},\int_{\mathbb{R}}{g_m
d\xi}\right]_{\alpha}-\left[\int_{\mathbb{R}}{f
d\xi},\int_{\mathbb{R}}{g_m d\xi}\right]_{\alpha}\right| & = &\dy
\left[\int_{\mathbb{R}}{f_n-f
d\xi},\int_{\mathbb{R}}{g_m d\xi}\right]_{\alpha},\\
& \leq & \dy \left\|\int_{\mathbb{R}}{f_n-f d\xi}\right\|_{\a} 
\left\|\int_{\mathbb{R}}{g_m d\xi}\right\|_{\a}^{\a-1},\\
& \leq & \dy \Psi_{\a}  \left\|\int_{\mathbb{R}}{g_m
d\xi}\right\|_{\a}^{\a-1} \int{|f_n-f|d\nu}.
\end{array}
$$
Since $(f_n)$ converges to $f$ in $L^1(\nu)$, we deduce that
$\dy \left[\int_{\mathbb{R}}{f_n d\xi},\int_{\mathbb{R}}{g_m
d\xi}\right]_{\alpha}$ converges to $\dy\left[\int_{\mathbb{R}}{f
d\xi},\int_{\mathbb{R}}{g_m d\xi}\right]_{\alpha}$. By the sequel
as $n$ tends to infinity in the equality (\ref{eq2-49})
we succeed to:
\begin{equation}\label{eq2-50}
\begin{array}{lll}
\displaystyle\left[\int_{\mathbb{R}}{f d\xi},\int_{\mathbb{R}}{g_m
d\xi}\right]_{\alpha} &=&\displaystyle
\sum_{j}\left(g^{m}_j\right)^{<\alpha-1>}
\int{f(\lambda)F_{B^{m}_j}(d\lambda)},\\
&=&\displaystyle \int{(g_{m}(\lambda'))^{<\alpha-1>}G(d\lambda')},
\end{array}
\end{equation}
where the measure $G$ is defined in the proposition
\ref{prop2-2}. Since $(g_m)_m$ converges to $g$ in $L_1(\nu)$, there exists a subsequence $(g_{m_k})_k$ of $(g_m)$ that
converges to $g$ ($\nu$-almost everywhere). By applying the fourth assertion of proposition \ref{prop2-2},
we have, $\left(g_{m}(\lambda')\right)^{<\alpha-1>}$ converges to
$(g(\lambda'))^{<\alpha-1>}$ ($v(G,.)$-almost everywhere).
Finally, since $G$ is of bounded total variations, then the Lebesgue convergence theorem ensure the convergence of  $\dy
\int{(g_{m}(\lambda'))^{<\alpha-1>}G(d\lambda')}$ to $\dy
\int{(g(\lambda'))^{<\alpha-1>}G(d\lambda')}$. By a similar technique as for $(f_n)$,
$$
\begin{array}{lll}
\dy\left| \left[\int_{\mathbb{R}}{f d\xi},\int_{\mathbb{R}}{g_m
d\xi}\right]_{\alpha}-\left[\int_{\mathbb{R}}{f
d\xi},\int_{\mathbb{R}}{g d\xi}\right]_{\alpha}\right| & \leq & \dy
\left\|\int_{\mathbb{R}}{f d\xi}\right\|_{\a} 
\left\|\int_{\mathbb{R}}{g_m - g d\xi}\right\|_{\a}^{\a-1},\\
& \leq & \dy (\Psi_{\alpha})^{\a-1}  \left\|\int_{\mathbb{R}}{f
d\xi}\right\|_{\a} \left(\int{|g_m-g|d\nu}\right)^{\a-1}.
\end{array}
$$
Since $(g_m)$ converges in $L^1(\nu)$ to $g$ then the
covariation, $\dy \left[\int_{\mathbb{R}}{f
d\xi},\int_{\mathbb{R}}{g_m d\xi}\right]_{\alpha}$ converges to
$\dy \left[\int_{\mathbb{R}}{f d\xi},\int_{\mathbb{R}}{g
d\xi}\right]_{\alpha}$. Finally, as $m$ tends to infinity
in the equality (\ref{eq2-50}) we obtain,
\begin{equation}
\displaystyle\left[\int_{\R}{f d\xi},\int_{\R}{g
d\xi}\right]_{\alpha}=\int{(g(\lambda'))^{<\alpha-1>}G_f(d\lambda')}
= I_1(f,g).
\end{equation}
The same reasoning on the integral  $I_2$ permits to show that $\dy \left[\int_{\mathbb{\mathbb{R}}}{f
d\xi},\int_{\mathbb{R}}{g d\xi}\right]_{\alpha} = I_2(f,g)$. In conclusion we have shown that the couple $(f,g)$ is MT-integrable and we have (\ref{eq2-48}).

{\bf Proof of proposition \ref{prop2-4}: }

Let $A$ and $B$ be any two Borel sets in $\R$ and suppose that there exist an other bimeasure $F'$ verifying (\ref{eq2-48}). In this case, by choosing the particular functions $f=\1_A$ and $g=\1_B$ we have :
$$
[d\xi(A),d\xi(B)]_{\a}=\int\int{\1_A(\lambda) \1_B(\lambda')
F(d\lambda,d\lambda')}=\int\int{\1_A(\lambda) \1_B(\lambda')
F'(d\lambda,d\lambda')}.
$$
which implies that $F(A,B)=F'(A,B)$ and the uniqueness of $F$ is proven.\\
Now let, $t_1,t_2, \dots, t_n$ finite sequence of reals. Let us show that the characteristic function of the S$\a$S random vector, $(X_{t_1}, \dots, X_{t_n})$ is characterized by $F$. For this purpose we distinguish the real and the complex cases.

\textbf{The real case: }  Using the definition of the characteristic function 
of real S$\a$S random vector, we have:
$$
\begin{array}{lll}
\dy\Phi_{X_{t_1},\dots,X_{t_n}}(\theta_1,\dots,\theta_n)& = &\dy \Phi_{\theta_1X_{t_1}+\dots+\theta_nX_{t_n}}(1),\\
& = &
\dy \exp\left\{-\left\|\theta_1X_{t_1}+\dots+\theta_nX_{t_n}\right\|_{\a}^{\a}\right\},\\
& =
&\dy \exp\left\{-\left[\theta_1X_{t_1}+\dots+\theta_nX_{t_n},\theta_1X_{t_1}+\dots+\theta_nX_{t_n}\right]_{\a}\right\}.
\end{array}
$$
Replacing $X_t$ by its expression given in (\ref{eq2-1}),
it is easy to see that,
$$
\begin{array}{lll}
\dy \left[\sum_{i=1}^{n}{\theta_i X_{t_i}},\sum_{i=1}^{n}{\theta_i
X_{t_i}}\right]_{\a}& = & \dy
\left[\sum_{i=1}^{n}{\theta_i\int{f(t_i,\lambda)d\xi(\lambda)}},
\sum_{i=1}^{n}{\theta_i\int{f(t_i,\lambda)d\xi(\lambda)}}\right]_{\a},\\
& = &
\dy\left[\int{\sum_{i=1}^{n}{\theta_if(t_i,\lambda)d\xi(\lambda)}},
\int{\sum_{i=1}^{n}{\theta_if(t_i,\lambda)}d\xi(\lambda)}\right]_{\a}.
\end{array}
$$
Since the functions $f(t_i,.)$ are bounded in
$\Lambda_{\a}(d\xi)$, which is the same as $\dy
\sum_{i=1}^{n}{\theta_if(t_i,.)}$, then according to (\ref{eq2-48}) of proposition \ref{prop2-3} we have:
$$
\left[\sum_{i=1}^{n}{\theta_i X_{t_i}},\sum_{i=1}^{n}{\theta_i
X_{t_i}}\right]_{\a}=\int\int{(\sum_{i=1}^{n}{\theta_if(t_i,\lambda)})
(\sum_{i=1}^{n}{\theta_if(t_i,\lambda')})^{<\a-1>}F(d\lambda,d\lambda')}.
$$
This implies that,
\begin{equation}\label{prob1}
\Phi_{X_{t_1},\dots,X_{t_n}}(\theta_1,\dots,\theta_n)=\exp\left\{-\int\int{(\sum_{i=1}^{n}{\theta_if(t_i,\lambda)})
(\sum_{i=1}^{n}{\theta_if(t_i,\lambda')})^{<\a-1>}F(d\lambda,d\lambda')}\right\}.
\end{equation}
It is clear from this last expression that the characteristic function is expressed directly by the known deterministic functions $f(t,.)$ and the unique bimeasure $F$.

\textbf{The isotropic complex case: }  It is known that isotropic complexes $\a$-stables random variables are parametric. Their characteristic function is entirely determined through the covariation norm $\|.\|_{\a}$, see for instance  \citet{Cam83, Uch99}.
Therefore, using the definition of the characteristic function of the random vector  $(X_{t_1},\dots,X_{t_n})$ we have:
$$
\begin{array}{lll}
\dy\Phi_{X_{t_1},\dots,X_{t_n}}(\theta_1,\dots,\theta_n)& = &\dy \Phi_{\theta_1X_{t_1}+\dots+\theta_nX_{t_n}}(1),\\
& = &
\dy \exp\left\{-c_0 \left\|\theta_1X_{t_1}+\dots+\theta_nX_{t_n}\right\|_{\a}^{\a}\right\},\\
& = &\dy \exp\left\{-c_0
\left[\theta_1X_{t_1}+\dots+\theta_nX_{t_n},\theta_1X_{t_1}+\dots+\theta_nX_{t_n}\right]_{\a}\right\},
\end{array}
$$
where $c_0$ is given by $\dy c_0=\frac{1}{2\pi}\int_0^{2\pi}{|\cos(\phi)|^{\a} d\phi}$, see [\citet{Sam94}, p.86]. Similarly as in the real case we have:
$$
\Phi_{X_{t_1},\dots,X_{t_n}}(\theta_1,\dots,\theta_n)=\exp\left\{-c_0  \int\int{(\sum_{i=1}^{n}{\theta_if(t_i,\lambda)})
(\sum_{i=1}^{n}{\theta_if(t_i,\lambda')})^{<\a-1>}F(d\lambda,d\lambda')}\right\}
$$
Which achieve the proof of the proposition.
\par

{\bf Proof of proposition \ref{prop2-5b} : } The proof of (2) of this proposition is the same as in \citet[proposition 1]{Hur89} for the case of periodically correlated processes. The proof of (3) is the same as in the almost periodically case see \citet[proposition 6]{Hur91}. In order to give the idea of the proof, we recall here the proof in the periodically covariated case. First, for all fixed $\t$, by using the representation  (\ref{eq2-48}) and the fact that $(e^{\imath t\l'})^{<\a-1>}=|e^{\imath t\l'}|^{\a-2}  e^{-\imath
t\l'} = e^{-\imath t\l'}$, we have:
\begin{displaymath}
\begin{array}{llll}
\displaystyle [X_{t+\tau}, X_{t}]_{\alpha}& = & \displaystyle
\int_{\mathbb{R}}\int_{\mathbb{R}}{e^{\imath (t+\tau) \lambda} 
e^{-\imath t \lambda'}F(d\lambda,d\lambda')},\\
& = & \displaystyle \int_{\mathbb{R}}\int_{\mathbb{R}}{e^{\imath
\tau \lambda}  e^{\imath
t(\lambda-\lambda')}F(d\lambda,d\lambda')}.
\end{array}
\end{displaymath}
On the other hand, according to the periodicity of the covariation function, for all $N>0$ we have :
\begin{displaymath}
\begin{array}{llll}
\displaystyle [X_{t+\tau}, X_{t}]_{\alpha}& = & \displaystyle
\frac{1}{2N+1}\sum_{-N}^{N}{[X_{t+\tau+kT}, X_{t+kT}]_{\alpha}},\\
& = & \displaystyle
\frac{1}{2N+1}\sum_{-N}^{N}{\int_{\mathbb{R}}\int_{\mathbb{R}}{e^{\imath
\tau \lambda}  e^{\imath t(\lambda-\lambda')} e^{\imath
kT(\lambda-\lambda')}
F(d\lambda,d\lambda')}},\\
& = & \displaystyle \int_{\mathbb{R}}\int_{\mathbb{R}}{e^{\imath
\tau \lambda}  e^{\imath
t(\lambda-\lambda')}D_N(\lambda,\lambda') F(d\lambda,d\lambda')}.
\end{array}
\end{displaymath}
with, $\dy  D_N(\lambda,\lambda')=
\frac{1}{2N+1}\sum_{-N}^{N}{e^{\imath kT(\lambda-\lambda')}}$. It is the  Fejer's kernel (see \citet{Hur91}). It is easy to see that:
\begin{equation}\label{eq2-54}
D_N(\lambda,\lambda')=\left\{
\begin{array}{llll}
 1 & \hbox{si} &  \lambda-\lambda'=\frac{2\pi k}{T}, \\
\frac{1}{2N+1}\frac{\sin((N+\frac{1}{2}) T (\lambda-\lambda'))}{\sin(T(\lambda-\lambda'))}
& \hbox{si}& \lambda-\lambda'\neq\frac{2\pi k}{T}.
\end{array}
\right.
\end{equation}
As $N$ tends to infinity in (\ref{eq2-54}),
$D_N(\lambda,\lambda')$ converge to $\1_L$ where $\displaystyle
L=\cup_{k\in\mathbb{Z}}{L_k}=\dy \bigcup_{k\in\mathbb{Z}}\left\{(\lambda,\lambda'),
\lambda-\lambda'=\frac{2\pi k}{T}\right\}$. Since  $F$ is of bounded Vitali variations, then by the Lebesgue convergence theorem which is true in our case, we have for all $t$ and $\tau$:
$$\int_{\mathbb{R}}\int_{\mathbb{R}}{e^{\imath (t+\tau) \lambda} 
e^{-\imath t
\lambda'}F(d\lambda,d\lambda')}=\int_{\mathbb{R}}\int_{\mathbb{R}}{e^{\imath
(t+\tau) \lambda}  e^{-\imath t
\lambda'}\1_{L}(\lambda,\lambda')F(d\lambda,d\lambda')}
$$
 Applying the inverse Fourier transform we see immediately that $F$ is concentrated on  
$L=\bigcup_{k\in\mathbb{Z}}{L_k}$. \\
Reciprocally, by using the fact that $F$ is concentrated on
$L$ we have:
\begin{displaymath}
\begin{array}{llll}
\displaystyle [X_{t+\tau+T}, X_{t+T}]_{\alpha}& = & \displaystyle
\int_{\mathbb{R}}\int_{\mathbb{R}}{e^{\imath (t+\tau) \lambda} 
e^{-\imath t \lambda'} e^{\imath T(\lambda-\lambda')}F(d\lambda,d\lambda')},\\
& = & \displaystyle \int_{\mathbb{R}}\int_{\mathbb{R}}{e^{\imath
(t+\tau) \lambda}  e^{-\imath t \lambda'}e^{\imath
T(\lambda-\lambda')}\1_L(\lambda,\lambda')F(d\lambda,d\lambda')}.
\end{array}
\end{displaymath}
However, for all ($\lambda,\lambda'$) in S, $T(\lambda-\lambda')$
belongs to $2\pi\mathbb{Z}$ then $e^{\imath
T(\lambda-\lambda')}=1$ by the sequel:
\begin{displaymath}
\displaystyle [X_{t+\tau+T}, X_{t+T}]_{\alpha} =  \displaystyle
\int_{\mathbb{R}}\int_{\mathbb{R}}{e^{\imath (t+\tau) \lambda} 
e^{-\imath t \lambda'} F(d\lambda,d\lambda')}=\displaystyle
[X_{t+\tau}, X_{t}]_{\alpha},
\end{displaymath}
which achieve the proof of the proposition.

\begin{proof}[ Proof of lemma \ref{lem2-1}:]
We have the following classic equality that one can find, for example, in \citet{Nik96}:
\begin{equation}\label{eq2-5}
\int_{0}^{+\infty}{\frac{1-\cos(s t)}{t^{\a+1}}dt}=|s|^{\a}
\frac{\Gamma(1-\a) \cos(\frac{\a\pi}{2})}{\a}.
\end{equation}
Since $1<\a<2$, for all real $s$, the map
$t\longmapsto\frac{\sin(st)}{t^{\a}}$ is integrable on
$[0,+\infty[$ because when $t$ tends toward infinity we have the asymptotic approximation 
$|\frac{\sin(t)}{t^{\a}}|\underset{\infty}{\sim} \frac{1}{t^{\a}}$
and when $t$ tends to $0$, $|\frac{\sin(t)}{t^{\a}}|
\underset{0}{\sim}\frac{1}{t^{\a-1}}$. We can then differentiate the equality (\ref{eq2-5}) under the integral sign with respect to $s$. We have therefore the following equality :
\begin{displaymath}
\int_{0}^{+\infty}{\frac{\sin(s t)}{t^{\a}}dt}=s^{<\a-1>}
\Gamma(1-\a) \cos(\frac{\a\pi}{2}).
\end{displaymath}
which achieves the proof of this lemma. {\qed}
\end{proof}

 \begin{proof}[Proof of lemma \ref{lem2-2}:]
 By using the complex notations $z=a+\imath b$ et $x=s+\imath t$ it is clear that
 $\RE(x \overline{z})=a s + b t$ and $|x|=\sqrt{s^2+t^2}$. Therefore, the calculation of the integral (\ref{eq2-6}) is equivalent to,
\begin{equation}\label{eq2-7}
I\triangleq \int_{-\infty}^{\infty}{\int_{-\infty}^{\infty}{\frac{1-\cos(a s
+ b t)}{(s^2+ t^2)^{\frac{p+2}{2}}}ds}dt}.
\end{equation}
To calculate the  integral I, for all $a$ and $b$, we use the change of variables:
\begin{equation}\label{eq2-8}
\left\{
\begin{array}{lll}
 \dy s = s(x,y) & = &\dy  \frac{1}{a^2+b^2}[(a+b) x + (a-b) y],\\
\dy  t = t(x,y) & = &\dy \frac{1}{a^2+b^2}[(b-a) x + (a+b)y],
\end{array}
\right.
\end{equation}
First, the Jacobian of the transformation (\ref{eq2-8}) is given by:
$$
J=\left| \begin{array}{lll}
 \dy \frac{a+b}{a^2+b^2}  &   & \dy \frac{a-b}{a^2+b^2}\\
\dy  \frac{b-a}{a^2+b^2}  &   & \dy \frac{a+b}{a^2+b^2}
\end{array}
\right| = \frac{2}{a^2+b^2}.
$$
A simple calculus shows that: $ a s + b t = x + y$ et $ s^2 + t^2
= \frac{2(x^2+ y^2)}{a^2+b^2}$, this implies:
\begin{equation}
\begin{array}{lll}
I & = & \dy
\int_{-\infty}^{\infty}{\int_{-\infty}^{\infty}{\frac{1-\cos(x+y)}{(\frac{2(x^2+
y^2)}{a^2+b^2})^{\frac{p+2}{2}}}  \frac{2}{a^2+b^2}dx}dy},\\
& = & \dy 2^{\frac{-p}{2}} (a^2+b^2)^{\frac{p}{2}}   
\int_{-\infty}^{\infty}{\int_{-\infty}^{\infty}{\frac{1-\cos(x+y)}{(x^2+
y^2)^{(\frac{p+2}{2})}} dx}dy}.
\end{array}
\label{eq2-9a}
\end{equation}
By a change of variables to the polar coordinates and thereafter using (\ref{eq2-5}), we have:
$$
\begin{array}{lll}
\dy
\int_{-\infty}^{\infty}{\int_{-\infty}^{\infty}{\frac{1-\cos(x+y)}{(x^2+
y^2)^{\frac{p+2}{2}}} dx}dy} & = & \dy \int_{0}^{2
\pi}{\int_{0}^{\infty}{\frac{1-\cos(r(\cos(\theta)+\sin(\theta)))}{r^{(p+2)}}
r dr}d\theta},\\
& = &\dy  \frac{\Gamma(1-p) \cos(\frac{p\pi}{2})}{p}  \dy \int_{0}^{2
\pi}{\left|\cos(\theta)+\sin(\theta)\right|^{2p}d\theta},\\
& = & \dy \frac{\Gamma(1-p) \cos(\frac{p\pi}{2})}{p} \int_{0}^{2
\pi}{\left|1+\sin(2\theta)\right|^{\frac{p}{2}}d\theta}.
\end{array}
 $$
By replacing this last equality in (\ref{eq2-9a}), we obtain (\ref{eq2-6}).
\par
Let us consider the derivative operator in space of complexes numbers $\mathbb{C}$ defined for  $z=a+\imath . b$ by, $\dy
\frac{\partial}{\partial z} = \left(\frac{\partial}{\partial a}-
\imath  \frac{\partial}{\partial b}\right)$. As in lemma \ref{lem2-1} we differentiate with respect to $z$ the equality
(\ref{eq2-6}) under the integral sign, we find:
\begin{equation}\label{eq2-11}
\int_{-\infty}^{\infty}{\int_{-\infty}^{\infty}{\frac{\partial}{\partial
z}\left(\frac{1-\cos(\RE(x \overline{z}))}{|x|^{(p+2)}}\right)ds}dt}=c(p)
 \frac{\partial}{\partial z} |z|^{p}.
\end{equation}
A simple calculus of the derivative show that, $\dy
\frac{\partial}{\partial z} |z|^{p} = p  z^{<p-1>}$ by using the same technique we have, $\dy \frac{\partial}{\partial z}
\cos(\RE(x \overline{z})) = -\overline{x}
\sin(\RE(x \overline{z})) $. We then use this results in 
(\ref{eq2-11}) we conclude that:
\begin{equation}\label{eq2-12}
\int_{-\infty}^{\infty}{\int_{-\infty}^{\infty}{\frac{\sin(\RE(x\overline{z}))}{(\overline{x})^{<p+1>}}ds}dt}=
p  c(p)  z^{<p-1>}.
\end{equation}
This achieve the proof of lemma \ref{lem2-2}.
\qed\end{proof}

\end{document}